# USING PIVOT SIGNS TO REACH AN INCLUSIVE DEFINITION OF RECTANGLES AND SQUARES

Maria G. Bartolini Bussi, Anna Baccaglini-Frank

Dipartimento di Educazione e Scienze Umane

Università di Modena e Reggio Emilia (Italia)

*We present some fragments of a teaching experiment realized in a first grade classroom, to sow the seeds for a mathematical definition of rectangles that includes squares. Within the paradigm of semiotic mediation, we studied the emergence of pivot signs, which were exploited by the teacher to pave the way towards an inclusive definition of rectangles and squares. This was done to favor overcoming children's spontaneous distinction of these figures into distinct categories, reinforced by everyday language. The experiment is an example of an approach towards the theoretical dimension of mathematics in early childhood.*

**Keywords:** first grade, semiotic mediation, pivot signs, inclusive definition, bee-bot

## INTRODUCTION

Rectangles and squares represent a paradigmatic example of the conflict between the perceptual experience and the theoretical needs of a mathematical definition (on this persisting conflict also see Hershkowitz 1990; Clements 2004; Fujita, 2012; Koleza & Giannisi 2013), where squares are to be considered as particular rectangles (we will refer to a definition of rectangles that includes squares as being *inclusive*). Mariotti and Fischbein (1997) claim that "from the figural point of view squares and non-square rectangles look so different that they impose the need of being distinguished at least as much as triangles and quadrilaterals" (Mariotti and Fischbein, 1997, p. 224). Actually the difficulty of naming and classifying geometrical figures (and, in particular, squares and rectangles), according to inclusive criteria, seems to depend on different reasons:

- the implicit constraints of everyday language: for instance, both in Italian and in English (as well as in other European languages) the names "quadrato" [square] and "rettangolo" [rectangle] hint at a complete separation of the figures into two different classes (square and not-square rectangles);
- some widespread improper practices in school which reinforce the separation between squares and rectangles (for instance, activities with attribute blocks, where squares and non-square rectangles are classified in different sets).

Hence, teaching needs to orient learning towards an inclusive definition. The question is: *at what age?* We claim that, although this choice may create a discontinuity between everyday language and school language, it is possible from early childhood to sow the seeds of an inclusive definition, focusing on the experience of walking along or drawing a rectangular path, where the change of direction in the four angle vertexes has the potential to attract the students' attention. In the following, we report on some fragments of a long term teaching experiment, carried out within the theoretical framework of semiotic mediation (Bartolini Bussi & Mariotti, 2008). Additional details are discussed by Bartolini Bussi and Baccaglini-Frank (in press).

**THEORETICAL FRAMEWORK**

In order to design and to analyze the teacher's role in the classroom teaching process, we adopted the theoretical framework of semiotic mediation (Bartolini Bussi & Mariotti, 2008; Bartolini Bussi, 2013). The design process is represented by the reciprocal relationships between the tasks, the artifact, and the mathematical knowledge at stake. In this relationship the *semiotic potential* of the artifact is made explicit. The artifact is the *bee-bot*, a small programmable robot represented in Figure 1 (also see the next section). When children are assigned a task they engage in a rich and complex semiotic activity, producing traces (gestures, drawings, oral descriptions and so on), that we refer to as "situated texts". The teacher's job is to collect all these traces (by observing and listening to the children), to analyze them and to organize a path for their development into "mathematical texts" that can be put in relationship with the fragments of mathematics knowledge that are to come into play. The process of semiotic mediation also concerns the functioning within the classroom. The teacher acts as a cultural mediator, in order to exploit, for all students, the semiotic potential of the artifact (the bee-bot in our case). In this last process, Bartolini Bussi and Mariotti (2008) identify three main categories of signs: artifact signs, pivot signs, and mathematical signs. *Artifact signs* "refer to the context of the use of the artifact, very often referring to one of its parts and/or to the action accomplished with it. […]"; *mathematics signs* "refer to the mathematics context" and *pivot signs,* which "refer to specific instrumented actions, but also to natural language, and to the mathematical domain" (ibid. p. 757). *Pivot signs* can be particularly useful for fostering a transition from situated "texts" to mathematical texts. Pivot signs develop and are enriched by their relationships with other pivot signs, hence building a *network* of pivot signs. Mathematical signs are not intended to suddenly substitute artifact signs; in fact the latter may survive for some time, especially for lower achievers or in cases in which the formal mathematical definition and the reasoning of the corresponding concepts require long term processes to be achieved.

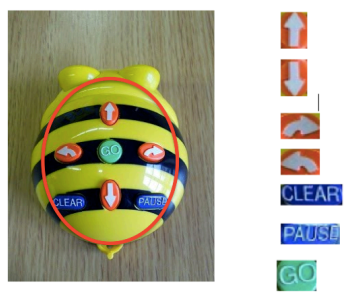

**Figure 1: Bee-bot's back and commands**

**THE CHOSEN ARTIFACT: THE BEE-BOT**

The bee-bot (Figure 1) is a small programmable robot, especially designed for young students. Its ancestor is the classical LOGO turtle, originally a robotic creature that could be programmed through an external computer to move around on the floor (LOGO Foundation, 2000). It is not necessary to have any additional computer to program the bee-bot; this can be done simply pressing a sequence of command buttons on its back. When the program is executed, the bee-bot moves on the floor: the execution of each command is followed by a blink of the eyes and by a short beep-sound. The bee-bot hints at many sets of meanings and mathematical processes, partly related to mathematics and partly related to computer science, for instance: counting (the commands); measuring (the length of the path, the distance); exploring space, constructing frames of reference, coordinating spatial perspectives, programming, planning and debugging. In a long term teaching experiment, all these sets of meanings are at stake, sometimes in the foreground and sometimes in the background. Focusing on any set of them depends on the adult's teaching intention. The bee-bot walks on the floor and traces paths that can be perceived, observed, described with words, gestures, drawings, sequences of command-icons and so on. Paths (either traced or imaginary, when no trace mark is actually left) constitute a large experiential base to "study" some plane figures, that can be traced using the available commands. These are polygons with sides measured by a whole number of steps and with right angles only. With the additional constraint of being convex, the bee-bot can be programmed only to turn "left" or "right" (with respect to itself), and therefore the convex polygons it can trace are always rectangles (including squares). Moreover, in experiences where "pretending to be the bee-bot" is essential, children embrace the robot's perspective: they move with the bee-bot and they see through its eyes. In particular, when walking along a closed convex path and ending up where they started, the children turn $360°$ in four equal "chunks" during which their orientation is perceived as essential (they find it important to end up facing the same direction as when they started).

**THE TEACHING EXPERIMENT**

Above we have discussed some features that define bee-bot's high semiotic potential with respect to the emergence of an inclusive definition of rectangles, characterized by the property of four right angles. Our teaching experiment was designed to capitalize on bee-bot's potential of fostering awareness of the "four right angles" property of generic rectangles (including squares).

Several sessions (15) were carried out in a first grade classroom at the beginning of the school year, for 4 months (more or less once a week) either in the classroom or in the gym, with a careful alternation of whole class or small group activity (with adult's guidance) and some individual activity. Each session was carefully observed by the teacher, by a student teacher or by a researcher (the second author of this paper), with the collection of students' protocols, photos, and videos. The tasks were designed by the whole research team (including all the mentioned adults), drawing on the initial intention and on some changes implemented "on the fly" based on episodes that occurred during the experiment. Due to space constraints it is not possible to report on all the details, so we have focused on particular sessions where the production of signs was very rich and fundamental for preparing the final summary texts and poster for the students (see Figure 6 in this paper, and Bartolini Bussi & Baccaglini-Frank, in press).

**Observing programmed bee-bots**

In this session, students were given two bee-bots that had ahead of time been programmed with the same sequence. The task was: *Describe what they do*. The students watched the twin bee-bots move together, starting facing in the same or in different directions, and then moving separately. Then the memory of one of the bee-bots was erased (CLEAR command-icon) and the students were asked to reprogram it so that it would move just like the other bee-bot. The students' productions concerned both global and local aspects. Global aspects refer to the perception of a path as a whole (as if bee-bot had drawn it on the floor), whilst local aspects refer to special points of the path. An example of the former is the expression "it did an L"; an example of the latter is "they switched the turn". Both aspects also appeared in gesturing: the path is represented by a single pointer finger tracing a path in the air (tracing gesture), whilst turning is represented by moving the right hand (for a right turn) or left hand (for a left turn) up and to the right or left in a rotation (turning gesture). The *turning gesture* was mirrored by the student-teacher, as a pivot sign with respect to the notion "angle" in a path.

**Pretending to be a bee-bot**

During this session the students were asked to work in pairs: one pretended to be the bee-bot and the other gave the first commands to move according to some undisclosed (to the first student) path. The intention was to guide the children to focus their attention on the turn command. Typical words used were be "Straight Ahead" "Left" "Right" "Backwards" usually without quantifying the number of steps, and frequently combining a translation with a change of direction (rotation). For example, when a student said "left" the bee-bot student frequently would not only turn left, but s/he would also take a step in that direction, or even just take a step to his/her left without even turning in that direction. The student-teacher's intervention here was fundamental in focusing the children's attention on "turn" commands, which led to their beginning to explicitly consider rotations as important elements per se, without having to associate them to steps.

**Constructing paths**

Several activities were designed around tracing different kind of paths on the floor. When the aim was to produce particular letters of the alphabet, the students' attention was focused mostly on the "possible" and "impossible" letters: they empirically discovered that some capital letters (e.g. L, T. I) could be traced out, whilst others could not (e.g. B, A, D, O). In fact, neither acute angles ("sharp points") nor circular arcs ("fat curves") could be traced by the bee-bot. Children produced many examples of combinations of words, gestures and drawings, aiming at distinguishing the shapes (letters) which could or could not be drawn. There was a particularly rich production of words such as "angles", "(fat) curves", "diagonals", "(sharp) tips/points", "broken lines" and of related gestures and drawings. Suddenly, within this experience, an important event took place; this will be the seed of an inclusive definition of rectangles.

**The main pivot sign: the "squarized" O**

In a small group the following exchange occurred:

> Student-teacher: …Did you do an O?
>
> Student: No. Then it could do like this this this and this [he gestures four consecutive right angles] a squarized O. Ah, then it can make a square!

We have translated a non-existing Italian word (*quadratizzato*) into a non-existing English word (*squarized*). Other students started talking about "squarized Os" and other possible "squarized letters", intending letters that include one or more squarized Os within them (e.g. P, B). These squarized 0s were acknowledged by the teacher and the research team as pivot signs, hinting at both the perceived path produced by the bee-bot (*artifact sign*) and at a *square* (a figure, interpreted as a *mathematical sign*). The importance of the four consecutive right angles suggested to orient children's attention towards this feature, that seemed to put in shade the length of each piece of the traced path (the sides) and to put in the foreground the four changes of direction, common to all squarized Os.

**Focusing on the four right angles**

In the students' complex experience, each right angle appeared with seemingly different meanings, that also affected the signs used. These, initially, were mainly dynamic and related either to the student pretending to be a bee-bot or to the bee-bot:

a) *Dynamic change of direction of the student pretending to be a bee-bot;*
b) *Dynamic change of direction of the bee-bot under the effect of the turn command.*

In both these cases, however, the angle was the external angle, i.e. the region swept by the gaze of either the student or the bee-bot while changing direction. When the researcher proposed to draw the paths in a "faster way: using a mark like the one the bee-bot would make if a marker were used", she chose to mirror a sign produced by a student "a turn like this" close to the turning point of the path (see Figure 2).

The sign had the potential to become a pivot sign with respect to the notion of "angle" (external angle): it recalls the command-icon on bee-bot's back, but it is somewhat decontextualized, since there appears to be no explicit mention to the bee-bot. In addition to these dynamic signs, as the teaching experiment went on, the children developed other signs, which lacked such dynamic components:

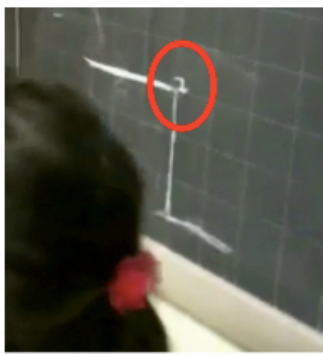

**Figure 2: Sign for the right angle**

c) *Hands-meeting gesture referring to the point in the path traced by the bee-bot;*
d) *Gestures to interpret a static figure (referring to a dynamic experience);*

e) Verbal utterance of the list of commands (uttered during or after the programming of the bee-bot);
  f) List of commands written horizontally.

First we describe the hands-meeting gesture (*type c*). While exploring figures that represented rectangles, including squares, a powerful gesture was realized by one of the groups of children and rapidly imitated by others: the two hands coming together at a right angle (Figure 3).

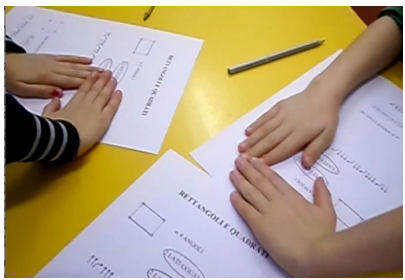

**Figure 3: the students' gesture**

The gesture emerged as the students tried to explain the property that all squarized O's (be they "allungati" [stretched] or "perfetti" [perfect]) had in common: all the four angles (internal angles) are equal and right. Moreover the gesture stresses the vertex as an important feature of the angle. Signs of *type d* were identified, for example, in the argument presented below (Figure 4a-e), on how the angles of a square or rectangle have to be (as opposed to angles such as the ones of the parallelogram that was included in one of the worksheets).

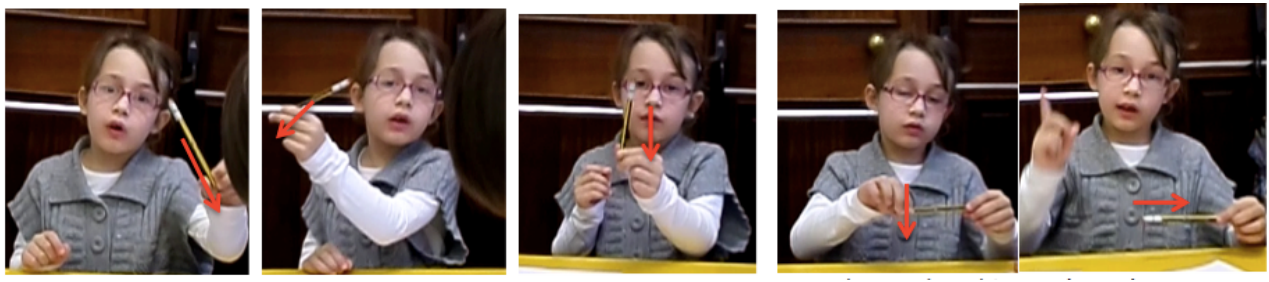

**Figures: 4a, 4b, 4c, 4d and 4e: Veronica's protocol**

Veronica is saying:

> [in a square or rectangle] the angles go down straight…[in the parallelogram] they are a bit down to the right and a bit down to the left. It has to go straight, not like this and down, it shouldn't be a bit down like this one [she moves her pencil in the air along a slanted line with respect to a horizontal bottom line, Figure 4a-c]. Instead it has to go straight like this and like this…it has to be straight like the line but a bit lying down [she marks the lower horizontal line, Figure 4d-e].

Signs of *type e* appeared when the students' attention was drawn to the "length of the path". Sometimes the turn command was in shade, as it did not lengthen the path perceived while the bee-bot spun around. However the number of commands for paths with angles, was not the same as the number of steps forward. So sometimes the turn command was still skipped (children 1, 2, 3, below). While in some of the children's utterances it was acknowledged (Child 4, below) as a command like the others (it is represented by a similar button and it is executed with by a beep and a blink of bee-bot's eyes).

Child 1: Three steps then three then three then three we make a square, because it is the same ends, the same length.

Child 2: Instead, the other one has 1,2-1,2,3-1,2-1,2,3, it has two the same and two the same.

Child 3: The other was three, two, three, two. Not all equal.

In contrast

Child 4: Two forward, turn right, two forward, turn right, two forward, turn right, two forward, turn right. The segments have to all be equal.

*Type f* emerged in activities in which children had learned to represent traced paths as written sequences of commands, typically in a horizontal line, from left to right. Within these sequences they searched for regularities allowing them to distinguish different types of "squarized Os". Figure 5 shows agreed-upon for the programmed sequences on the interactive white board after a discussion on "stretched squarized Os" (non-square rectangles) with respect to "perfect squarized Os" (squares).

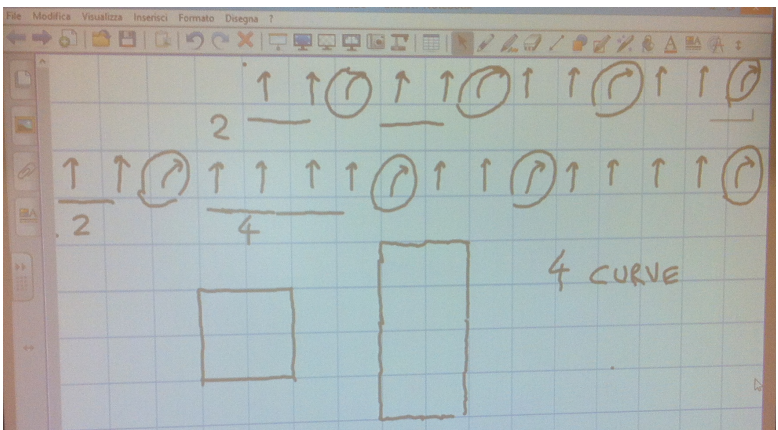

**Figure 5: Agreed-upon signs**

We note here how some students' language (in this and other occasions) seemed to be evolving into condensed pre-algebraic forms, such as a+b+a+b, that could eventually become expressions like 2a+2b for the rectangle and 4a for the square (a particular case in which a=b). In this teaching experiment, however, we did not pick up on these expressions, leaving them only as little germs to be nurtured by the teacher in future years (perhaps even during the second grade).

**Focus on the shapes as wholes**

Shapes as wholes were focused on from the very beginning of the teaching experiment, with either verbal descriptions alone or also with hand gestures. After the introduction of the idea of squarized Os, the adults involved in the experiment started mirroring students' utterances involving the words "rectangles" and "squares". As expected, when the attention was not brought to the word *squarized O* students spontaneously tended to partition the two situations, implying that "rectangles" had pairs of sides with *different* lengths ("equal in front of each other") while "squares" had sides that were "all equal". For some children this property seemed to persist when talking about "stretched squarized Os" with respect to "perfect squarized Os", while other children seemed to only differentiate "perfect squarized Os" from all other squarized Os, since they were special, being "all equal".

## The shared meanings

We chose to build on what seemed to be the idea of this second group of students to reach a summary of the shared meanings. The most important step in this direction was a poster of "our" discoveries, a first step towards the development of "mathematical texts".

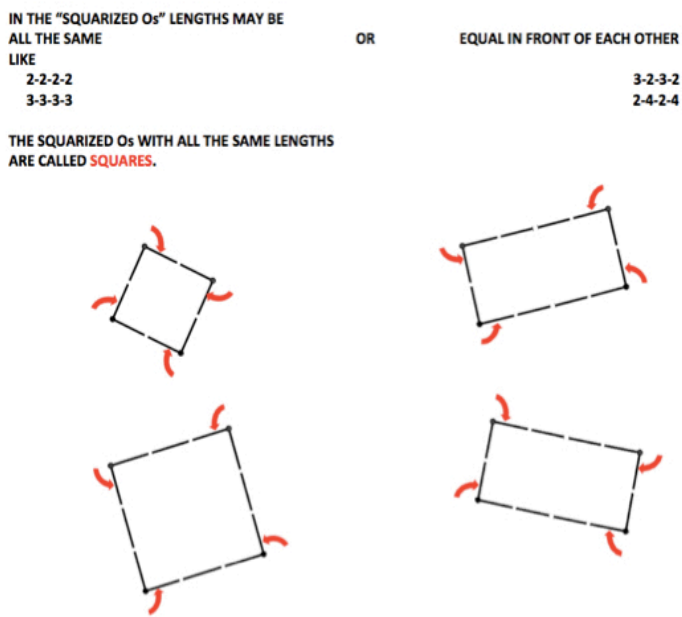

**Figure 6: Poster of "our" discoveries**

In this poster (Figure 6) several signs produced in the classroom are reconsidered, constructing a text where artifact signs (e.g. the figure of the bee-bot, the recollection of "giving it" a sequence of commands, the turns), pivot signs (e.g. the squarized Os, the small arrow to represent the external angles, and mathematical signs (e.g. squares; numbers; rectangles) are included.

Is this text a mathematical text? Not yet: it is still an *hybrid text*, where the richness of the exploration remains present. What is important in this phase is that *all* of the students could identify this poster as having been produced by the whole class as a community. The choice of which signs to include was discussed by the research team, trying to collect signs that hinted at the individual and collective processes. The poster was discussed in the classroom; the students seemed very happy to find their ideas made public and to receive a reduced-size copy to glue on their notebooks. Some months later, a follow up questionnaire confirmed that (at least some) students had appropriated, and transferred it to a mathematical context, an inclusive definition of rectangles (other students were still on their way along this process). As mentioned before, the process is not to be considered finished. The teacher has planned to go on with the same group of students and deepen the inclusive definition for which she planted the seeds during this teaching experiment in the first grade.

## DISCUSSION

The teaching experiment fruitfully exploited the semiotic potential of the bee-bot, joining different ways of representing the paths traced by the small robot, as sequences of commands, as wholes, as either physical or mental drawings, in both dynamic and static ways. During this long term process the students approached several pieces of mathematics knowledge, including counting (the commands), measuring (the length of the path, the distance), exploring space, constructing and changing frames of reference, coordinating spatial perspectives, programming, planning and debugging. The approach towards an inclusive definition of rectangles is only one aspect of this long and complex process.

The observation of the process has shown a very rich production of pivot signs, in the form of gestures, drawings, and utterances. Their shared feature is the manifold hint towards specific instrumented actions, natural language and mathematics language. These are fundamental elements of the students' process of meaning-making, hence they are exploited by the adults (e.g. mirrored) to foster evolution of the process itself. The special focus on gestures is consistent with Arzarello et al.'s (2008) theoretical construct of *semiotic bundle*, that allows focusing on the relationships of gestures with the other semiotic resources within a multimodal approach, and that helps to frame the mediating action of the teacher in the classroom.

A final comment on language. We do not claim that the inclusive (and decontextualized) definition of rectangles is already accepted by all the students (in fact we saw that this was not the case). Rather we find it important that students started becoming aware of the fact that theoretical mathematical needs may be different from everyday life needs. Moreover, we do not believe that the inclusive definition should be used also in everyday life. Rather it seems that, with this experiment, we have put the students in the situation of potentially seeing squares and rectangles within a same "family". What happened, indeed, was that the idea of "square" seemed to be overarching, in spite of the mathematical choices. The students seem to speak of the squarized O as the ancestor of rectangles (including squares) but, from the perceptual point of view they need to distinguish "perfect squares" from "stretched squares". This reminds us of the Chinese way of naming squares and rectangles: the sequences of ideograms for the words "square" and rectangle" contain two out of three of the same ideograms. Those that indicate "sides" and "shape" are the same, while the first indicates "exact" (for the square) and "long" (for the rectangle). This is represented in Figure 7.

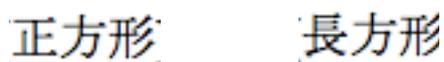

**Figure 7: Chinese characters**

**for "square" and "rectangle".**

So, linguistically, a square is seen as a "shape with exact sides" and a rectangle as a "(same) shape with long sides". In this case language makes explicit that squares and rectangles are two kinds of a same thing, deeply related to each other and not partitioned into categories. The Chinese choice of the square as the most important shape may be related to the Chinese ancient culture, where it represents the Earth and the circle represents the Sky. This iconic cosmology is shared by other ancient cultures.

ACKNOWLEDGEMENTS

This teaching experiment was carried out within the PerContare project (ASPHI, 2011), with the support of the Fondazione per la Scuola of the Compagnia di San Paolo di Torino.